\input amstex

\magnification=\magstep1
\input amsppt.sty
\hsize=14.4 cm
\hoffset -0.5cm
\overfullrule =10pt

\topmatter
\title Kobayashi hyperbolicity of almost complex manifolds\endtitle
\author Debalme R.\endauthor 
\abstract Our goal here is to give a simple proof of the non-integrable version of 
Brody's characterisation theorem. A different proof is given in  [KrOv].
\endabstract
\footnotetext {AMS subject classification 32 H 20, 53 C 15 \hfill \break
address: UFR de Math\'ematiques, Universit\'e Lille 1, 59655 Villeneuve d'Ascq Cedex, France.
\hfill \break e-mail: debalme\@gat.univ-lille1.fr}

\endtopmatter

\subhead 0. Introduction \endsubhead

In this note our aim is to extend the notion of Kobayashi hyperbolicity to 
the case of the almost complex manifolds.

Recall that an almost complex manifold is a pair $(M,J)$ where $M$ is a real manifold, 
countable at infinity and $J$ is a $C^1$-smooth section of $End(TM)$ with $J^2 = -Id$.
Denote by $J_{st}$ a standart almost complex structure in $\Bbb R^{2n}, n \geq 1$. 
Thus $(\Delta, J_{st})$ is a unit disk and $(\Bbb R^2,J_{st})$ is the usual complex line.

Let $(S,J_S)$ be a Riemann surface. Recall that a $C^1$-smooth mapping $u : (S,J_S)
\to (M,J)$ is called $(J_S,J)$-holomorphic if its differential commutes with the almost complex 
structures, i.e. $du \circ J_S = J \circ du$ as mappings from $TS$ to $TM$ (Cauchy-Riemann 
equation). If the structure 
$J_S$ is clear from the context we shall simply call $u$ $J$-holomorphic.

We shall prove, see Lemma 1, that provided $J \in C^2$, for any two points sufficently close 
$p,q \in (M,J)$ 
there exists a $J$-holomorphic map $u : \Delta \to (M,J)$ with $u(0)=p$ and 
$u({1 \over 2})= q$. This enables us to define the Kobayashi 
pseudo-distance $k^J_M$ on $(M,J)$. $(M,J)$ is said to be Kobayashi hyperbolic if $k^J_M$ is 
actually a distance.

Our goal is to prove the Brody characterisation theorem in the almost complex context.

\proclaim {Theorem} A compact almost complex manifold $(M,J)$, with $J \in C^2$, 
 is not hyperbolic if and only if it contains a non trivial $J$-complex line 
$f : \Bbb C \to (M,J)$.
\endproclaim

A line $f : \Bbb C \to (M,J)$ is not trivial if $f \not \equiv const$.

\subhead 1. Local properties and definition of Kobayashi pseudo-distance \endsubhead

\smallskip

It is known, see for example [Sk], theorem 3.1.1 that if $v$ is a small enough vector in $T_pM$,
then there exists a $J$-holomorphic map $ u : \Delta \to M$ with $u(0)=p$ and 
$du({\partial \over \partial z})=v$. We shall need a slightly modified version of this
 statement .

\proclaim {Lemma 1}  Let $(M,J)$ be an almost complex manifold, $J \in C^2$.
 If $p$ and $q$ are two
 points of $M$ sufficently near, there exists a $J$-holomorphic curve $ u : \Delta \to M$
such that $p$ and $q$ are in $u(\Delta)$.\endproclaim

\demo{ Proof} The issue being local, via a chart 
we can take $M=\Bbb R ^{2n}$. We can assume that $J (0) = i Id$.
 We are working in a ball 
$B=B(0,r) \subset \Bbb R ^{2n}$ sufficently small to have $J (v) + i$ invertible.
The equation expressing that $u$ is $J$-holomorphic is 
$$ {\partial u \over \partial y} =  J (u) {\partial u \over \partial x}.$$
Using the identities \quad
${\partial u \over \partial y} = {i \over 2}({\partial u \over \partial z} -
{\partial u \over \partial \bar z})$ \quad
and \quad ${\partial u \over \partial x} = {1 \over 2}({\partial u \over \partial z} +
{\partial u \over \partial \bar z})$ \hfill \break
one can rewrite the equation in the following form

$$ (i + J(u)){\partial u \over \partial \bar z} = (i- J(u)){\partial u \over \partial z}. $$
Since $i +J(u)$ is invertible, this can be written as

$$ {\partial u \over \partial \bar z} = q_{J}(u){\partial u \over \partial z} $$
Here $q_{J} : B \longrightarrow End_{\Bbb R}(\Bbb R ^{2n})$ is defined by
$ q_{J} (v) = [i +J (v)]^{-1} [i -J (v)]$.
We introduce the mapping from $L^{2,p}(\Delta,\Bbb R ^{2n})$ to
$L^{3,p}(\Delta,\Bbb R ^{2n}), p>2 :$
$$ P : \varphi \longmapsto (P\varphi)(z)={1 \over 2i\pi} \int _\Delta {\varphi(\zeta) \over
\zeta - z} d\zeta \wedge d\bar{\zeta} $$
the Cauchy-Green operator which verifies ${\partial \over \partial \bar{z}} \circ P = Id$.
Let us consider the map 

$$\Phi : (\varepsilon , u) \longmapsto [Id - Pq_{J}(\varepsilon u)
{\partial \over \partial z}] u $$

from $[0,1] \times L^{2,p}(\Delta,B)$ to $L^{2,p}(\Delta,\Bbb R ^{2n})$.
\smallskip

\noindent One has 
$${\partial \Phi (\varepsilon ,u) \over \partial \bar z} = {\partial u\over \partial \bar z}
- q_{J}(\varepsilon u) {\partial u\over \partial z}. $$
We see that $\Phi$ is of class $C^1$ 
and that $\Phi (\varepsilon ,u)$ is holomorphic in the standard sense if and only if 
$\varepsilon u$ is $J$-holomorphic. Denote $\Phi_{\varepsilon}=\Phi(\varepsilon,.)$. 
Since $q_{J}(0)=0$, $\Phi_0 = Id_{L^{2,p}(\Delta ,B)}$.
Thus $\exists \varepsilon _0 > 0 $ such that $ \forall \varepsilon \in [0,\varepsilon _0[,
\quad \Phi_{\varepsilon}$ is a diffeomorphism from $W$ to $V$ neighbourhoods of $0$ in 
$L^{2,p}(B)$ and $ L^{2,p}(\Bbb R ^{2n})$ respectively.
Consider the mapping from $\Delta$ to $\Bbb R ^{2n}$
$$ h_{p,q}  : \quad  z  \longmapsto p + 2z(q-p)  $$
where $(p,q) \in (\Bbb R ^{2n})^2$ and denote $u_{\varepsilon ,p,q} 
= \Phi_\varepsilon^{-1} h_{p,q}$. 
We remark that  \hfill \break
$\cdot$ $h_{p,q}$ being holomorphic, $\varepsilon u_{\varepsilon ,p,q} $
is $J$-holomorphic. \hfill \break
$\cdot$ $u_{0,p,q} = h_{p,q}$ . So it verifies $u_{0,p,q}(0)=p$ and
 $u_{0,p,q}({1 \over 2})=q$. \hfill \break
Consider the mapping from $[0, \varepsilon] \times (\Bbb R^{2n})^2$ to $(\Bbb R^{2n})^2$
 $$ \Xi : \quad (\varepsilon,p,q) \longmapsto (u_{\varepsilon,p,q}(0),u_{\varepsilon,p,q}
({1 \over 2}))$$
$\Xi$ is $C^1$ and from our last remark $\Xi(0,.,.) = Id_{(\Bbb R^{2n})^2}$.
So by the implicit functions theorem, if $\varepsilon$ is sufficently small
, there exist $U$ and $U^{'}$ neighbourhoods of zero in $(\Bbb R^{2n})^2$ such that
$ \Xi (\varepsilon,.,.) :  U  \longrightarrow U^{'}$
is a diffeomorphism.
Let $p_0$ and $q_0$ be two points sufficently near of zero (i.e. 
$({p_0 \over \varepsilon},{q_0 \over \varepsilon}) \in U^{'}$). 
There exists $(p,q)$ such that $ \varepsilon u_{\varepsilon,p,q}(0)=p_0$ and 
$\varepsilon u_{\varepsilon,p,q} ({1 \over 2}) = q_0$ . We have thus made 
$ \varepsilon u_{\varepsilon,p,q}$  a $J$-holomorphic curve which is going 
through $p_0$ and $q_0$.
\smallskip
\hskip 12 cm $\square$
\enddemo

\smallskip

This property of $J$-holomorphic curves enables us to define the Kobayashi pseudo-distance
on the almost complex manifolds.
Let $(M,J)$ denote an almost complex manifold, $J \in C^2$, and let $\rho$ 
be the Lobatchevski distance on $\Delta$. The associated metric is
$$ \rho = {dz \otimes d\bar z \over (1-|z|^2)^2}.$$

We define the Kobayashi pseudo-distance $k_M^J$ on $M$ as follows.
Given two points $p$ and $q$ in $M$, we choose points $p=p_0, p_1, \ldots ,p_{k-1}, p_k=q$
of $M$ and points $a_1, \ldots ,a_k, b_1, \ldots , b_k$ of $\Delta$ and $J$-holomorphic
mappings $f_1, \ldots , f_k$ from $\Delta$ to $M$ such that $f_i(a_i) = p_{i-1}$ and 
$f_i(b_i)=p_i$ for $i=1, \ldots , k$. For each choice of points and mappings made this way,
 we consider the number $\sum_{i=1}^k \rho (a_i,b_i)$.  
 Let $k_M^J(p,q)$  be the infinimum of the numbers obtained this way for all possible choices.
It is an easy matter to verify that $k_M^J : M \times M \longrightarrow \Bbb R$ is continuous
 and satisfies the axioms of the pseudo-distances : $k_M^J(p,q) \geq 0 ,
\quad k_M^J(p,q)=k_M^J(q,p) , 
\quad k_M^J(p,q) + k_M^J(q,r) \geq k_M^J(p,r) $.\hfill 

\smallskip
The following property of $k_M^J$ is obvious :
\proclaim{ Property 1} Let $f : (M,J) \longrightarrow (N,J^{'})$ be a 
$(J,J^{'})$-holomorphic mapping. Then $\forall (p,q) \in M^2$, one has
 $ k_M^J(p,q) \geq k_N^{J'}(f(p),f(q))$
\endproclaim

\proclaim{ Corrolary} $k_{\Bbb C} \equiv 0$
\endproclaim

\definition{ Definition } An almost complex manifold $M$ is said to be hyperbolic if $k_M^J$
 is actually a metric.
\enddefinition

\smallskip
We shall need the following 
\smallskip

\proclaim{Lemma 2} Let $(M,J)$ be an almost complex manifold.
The set of all $J$-holomorphic mappings $\Delta \to M$ is closed in the compact-open
topology.\endproclaim

\demo {Proof} Let $(f_n)_{n\in \Bbb N}$ a sequence of $J$-holomorphic mappings from 
$\Delta$ to $(M,J)$, converging uniformly on each compact of $\Delta$.
Denote $f$ the limit. 
Choose two compacts $K$ and $K^{'}$ of $\Delta$ such that $K$ is included in the interior of 
$K^{'}$.
From Sikorav [Sk] proposition 2.3.6 (i), p.171, we see that if $K^{'}$ is sufficiently small then
$\Vert f_n \Vert_{C^2(K)} \leq  L \Vert f_n \Vert_{L^{\infty}(K^{'})}$.
So uniform convergency of $(f_n)_{n\in \Bbb N}$ implies $C^2$-convergency to $f$.
 Thus  $f$ verifies the Cauchy-Riemann equations and is so $J$-holomorphic.
\smallskip
\hskip 12 cm $\square$
\enddemo

\subhead 2. Characterization theorem \endsubhead
\smallskip
\noindent
We suppose $M$ compact in all the following.
We denote by $\Cal O ((M,J);(N,J^{'}))$ the set of all $(J,J ^{'})$-holomorphic 
mappings from $M$ to $N$.
Fix some Riemannian metric $|.|$ on $M$. We start with the following

\proclaim{ Lemma 3}  $M$ is hyperbolic iff
$\sup \{ |f^{'}(0)| ,  f \in \Cal O (\Delta ; (M,J))\} < \infty $ 
where $f^{'}(z_0)=
f_{\ast}(z_0).{\partial \over \partial x}$\endproclaim

\demo{ Proof :  Sufficent condition} Observe that 
$$|f^{'}(0)|= |f_{\ast}(0).{\partial \over \partial x} |=|f_{\ast}(0).\nu|$$
We will note it $|f_{\ast}(\nu)|$. 
In the  precedent expression $\nu \in T_0 \Delta, \rho (\nu)=1$ is arbitrary. 
Indeed, at zero, $ \rho$ is the euclidian metric, so $\rho (\nu) = 1 $ means that 
$\nu = a {\partial \over \partial x} $ with $a$ a rotation. Thus
$$ \sup \{|f^{'}(0)|, f \in \Cal O(\Delta ,(M,J))\} =  
 \sup \{|f_{\ast}(\nu)|, f \in \Cal O(\Delta ,(M,J)); \nu \in T_0M, \rho (\nu)=1\} 
\eqno {(1)}$$
Let $\mu \in T_p\Delta$ , and let $\varphi_{0,p}$ be the conform automorphism of $\Delta$ 
exchanging $0$ and $p$. There exists $\nu \in T_0M$ such that $\varphi_{\ast}(\nu)= \mu$ and so
$f_{\ast}(\mu) = (f \circ \varphi)_{\ast}(\nu).$ \hfill \break
So
$$ \sup \{|f^{'}(0)|, f \in \Cal O(\Delta ,(M, J))\} =  
 \sup \{|f_{\ast}(\nu)|, f \in \Cal O(\Delta ,(M,J)); \nu \in TM, \rho (\nu)=1\}$$ 
We assume that $\sup \{|f^{'}(0)|, f \in \Cal O(\Delta ,(M,J))\} = c < \infty  $, 
and we want to prove that $k_M^J(p,q) . c \geq |p,q|$. Then we would be able to conclude 
$k_M^J(p,q)=0 \Leftrightarrow p=q $ (so $M$ is hyperbolic).
$$\eqalign {
|p,q| & = \inf \{ \int_{\gamma} |\gamma^{'}(s)|ds, \quad \gamma : [0,1] \to M , \gamma (0)=p
, \gamma (1)=q \} \leq \cr
& \leq \inf \{ \int_{\gamma} |\gamma^{'}(s)|ds, \quad \gamma : [0,1] \to M , \gamma (0)=p
, \gamma (1)=q, \gamma = \sum_{i=1}^k f_i(\delta_i) \cr
& f_i \in \Cal O(\Delta;(M,J)), \delta_i: [0,1] \to \Delta \} \leq \cr
& \leq \inf \{ \sum_{i=1}^k \int_0^1 f_{i \ast}(\delta_i^{'}(s) )ds
\quad f_i \in \Cal O(\Delta;(M,J)), f_i(\delta_i(1))=f_{i+1}(\delta_{i+1}(0)),
 f_1(\delta_1(0))=p,\cr
& f_k(\delta_k(1))=q \} \leq \cr 
& \leq \inf \{ \sum_{i=1}^k c \rho (\delta_i),
\quad f_i \in \Cal O(\Delta;(M,J)), f_i(\delta_i(1))=
f_{i+1}(\delta_{i+1}(0)),f_1(\delta_1(0))=p, \cr
& f_k(\delta_k(1))=q \} \cr
& = c k_M^J(p,q) \cr } $$
 
Here in the first inequality, we shifted to the curves piecewise lying in $J$-holomorphic
 disks. In the third, we used the fact that $|{f_i}_{\ast}(\delta_i^{'}(s))| \leq
 C \rho(\delta_i^{'}(s))$ by (1).

\smallskip

 \noindent {\it Necessary condition}.
Let us assume that $\sup \{|f^{'}(0)|; f \in \Cal O(\Delta;(M,J)) \} = \infty$. \hfill \break 
Then there exists a sequence of $J$-holomorphic mappings $(f_n)_{n\in \Bbb N}$, such that
$ \displaystyle\lim_{n \to \infty}|f_n^{'}(0)| = \infty$.
M being compact, we extract from this sequence a converging subsequence : 
$\displaystyle\lim_{n \to \infty} f_n(0) = p$. 
Let $U$ be a coordinate neighbourhood of $p$. We suppose that $J(p)=i.Id$ and we take $U$ small
enough to have $ \Vert q_J(x) \Vert < \varepsilon$ for $x \in U$, where $ \varepsilon > 0$
comes from proposition 2.3.6 (i) of [Sk], p.171.
From this proposition of Sikorav, we know that for any
$f \in \Cal O(\Delta;(U,J))$ 
$\Vert f \Vert_{C^2(\Delta_{{1 \over 2m}})}
\leq  L_m \Vert f \Vert_{L^{\infty}(\Delta_{{1 \over m}})}$  
if $f(\Delta_{1 \over m}) \subset U$,   
where $L_m$ only depends on $m$.
As $\Vert f \Vert_{C^2(\Delta_{{1 \over 2m}})}=\Vert f \Vert_{L^{\infty}(\Delta_{{1 \over 2m}})}
+\Vert df \Vert_{L^{\infty}(\Delta_{{1 \over 2m}})}
+\Vert d^2f \Vert_{L^{\infty}(\Delta_{{1 \over 2m}})}$
one has if $f_n(\Delta_{{1 \over m}}) \subset U$

$$|f^{'}_n(0)|\leq L_m \Vert f_n \Vert_{L^{\infty}(\Delta_{{1 \over m}})}$$

So, for any $m$ there exists $n_m$ such that 
$f_{n_m}(\Delta_{{1 \over m}}) \cap \partial U \not= \varnothing $.
In particular we may choose a sequence of points $(x_m)_{m \in \Bbb N}, x_m \in
 f_{n_m}(\Delta_{{1 \over m}})$ and therefore $k_M^J(f_{n_m}(0),x_m) \leq 
\rho(0,{1 \over m})$ converges to zero.
Since $k_M^J$ is continuous and $\partial U$ compact $(x_m)_{m \in \Bbb N}$ has an accumulation
point $x$ which verifies $k_M^J(x,p)=0$, with $x \not= p$ and $x$ in $\partial U$.
So $k_M^J$ is not a metric.
\smallskip
\hskip 12 cm $\square$

\enddemo

We remark now that the reparametrisation lemma of Brody (cf [Br]) remains valid in the
 nonintegrable case.

\proclaim{ Lemma 4} Let $(M,J)$ be an almost complex manifold. 
Let $f : \Delta_r \to M$ be a $J$-holomorphic curve  with $|f^{'}(0)| \geq c \geq 0$.
Then there exists $\tilde f$ a $J$-holomorphic curve $\tilde f :\Delta_r \to M$ such that
$$\displaystyle\sup_{z \in \Delta_r} |\tilde f^{'}(z)|({r^2-|z|^2 \over r^2}) =
|\tilde f^{'}(0)| = c$$ \endproclaim 

\demo {Proof} First we will arrange to have equality, and then we will force the supremum to 
occur at the origin. For $t \in [0,1]$ let $f_t : \Delta_r \to M $ be the map 
$z \mapsto f(tz) $.
Let $s(t) = \displaystyle\sup_{z \in \Delta_r} |f^{'}_t(z)|({r^2-|z|^2 \over r^2})$.
Then for any $t < 1, \quad \displaystyle\sup_{z \in \Delta_r} |f^{'}_t(z)| \leq 
\displaystyle\sup_{z \in t\Delta_r} |f^{'}(z)| \leq
\displaystyle\sup_{z \in \overline{\Delta_{tr}}} |f^{'}(z)| < \infty$.
Since $f$ is continuous on $\overline{\Delta_{tr}}$, and $({r^2-|z|^2 \over r^2}) \leq 1$,
 we have $s(t) < \infty $ for $ t < 1$. \hfill \break
Since $t \mapsto \displaystyle\sup_{z \in \overline{\Delta_{tr}}} |f^{'}(z)| $ is continuous
$s$ is also continuous. 
From $s(0)=0$ and $\displaystyle\lim_{t \to 1}s(t) \geq c $ 
we deduce that there exists
$t_0$ in $[0,1]$ such that $s(t_0)=c$.

\medskip
\noindent {\it First case : $t_0 = 1$.} 
$$c=s(1)= \displaystyle\sup_{z \in \Delta_r} |f^{'}(z)|({r^2-|z|^2 \over r^2})$$
While for $z=0 \quad |f^{'}(0)|{r^2 \over r^2} \geq c $ \quad the supremum occurs at $z=0$.
Just take $\tilde f = f$.

\medskip
\noindent {\it Second case : $t_0 < 1$.}\quad
The supremum is reached at a point $z_0$ inside $\Delta_r$. Let L be the conform automorphism
of $\Delta_r$ exchanging $0$ and $z_0$. Denote $\tilde f = f_{t_0} \circ L$.
Since  the quantity $|f^{'}_{t_0}(z)|({r^2-|z|^2 \over r^2})$ measures the derivative with
respect to $\rho_r$, it is invariant under $L$.
\smallskip
\hskip 12 cm $\square$
\enddemo

\proclaim {Lemma 5} Let $M$ an compact manifold. The familly $\Omega$ of all 
$C^{\infty}(\Delta_r,M)$ mapping satisfying 
$$\displaystyle\sup_{z \in \Delta_r} |f^{'}(z)|({r^2-|z|^2 \over r^2}) = |f^{'}(0)| = c$$
is relatively compact in the compact-open topology \endproclaim

\demo{ Proof} This familly is equicontinuous. 
Let $(z_n)_{n\in \Bbb N}$ be a dense sequence of points of $\Delta$.
Let $(f_n)_{n\in \Bbb N}$ be a sequence of $\Omega$.
 \quad $(f_m(z_1))_{m\in \Bbb N}$ is a sequence in $M$ which is compact. 
A converging subsequence can be extracted. 
From this subsequence we can extract another, which converges at $z_2$.
We make that way a subsequence of $(f_n)_{n\in\Bbb N}$ which converges for all $z_n,
 n \in \Bbb N$(diagonal process).
Since $(z_n)_{n\in \Bbb N}$ is dense and $(f_n)_{n\in \Bbb N}$ is equicontinuous, 
$(f_n(z))_{n\in \Bbb N}$ converges for all $z$ in $\Delta$.
Equicontinuity tells us that $(f_n)_{n\in \Bbb N}$ converges uniformly over compact subsets 
to this limit, which is so continuous.
\smallskip
\hskip 12 cm $\square$
\enddemo

\demo {Proof of the characterization theorem. \quad Sufficent condition}
Assume $M$ contains a non trivial almost complex line $f \in \Cal O(\Bbb C;(M,J))$.
Consider two distinct points $p$ and $q$ of $M$ such that $p=f(x)$ and $q=f(y)$.
$$k_M^J(p,q) \leq k_{\Bbb C}(x,y)$$
Since $k_{\Bbb C} \equiv 0, \quad k_M^J(p,q)=0$ and $k_M^J$ is not an actual distance.

\medskip
\noindent {\it Necessary condition.}
Assume $M$ is not hyperbolic. Lemma 3 tells us that there exists $(f_n)_{n\in \Bbb N}$ a 
sequence of $ \Cal O(\Delta;(M,J))$ such that $\displaystyle\lim_{n \to \infty}
|f_n^{'}(0)|= \infty$. 

Consider the mappings $g_n : \Delta_{r_n} \to M $ given by $z \longmapsto f_n({z \over r_n})$
where $r_n =|f_n^{'}(0)|$.
We have $|g_n^{'}(0)|=1 $.
From lemma 4, we obtain a sequence $(\tilde g_n)_{n \in \Bbb N}$ of $J$-holomorphic curves 
$\Delta_r \to M$ such that $\displaystyle\sup_{z \in \Delta_r} 
|\tilde g_n^{'}(z)|({r^2-|z|^2 \over r^2}) = |\tilde g_n^{'}(0)| = 1$.
From lemma 5 we can extract a subsequence of $(\tilde g_n)_{n \in \Bbb N}$ which 
converges on $\Bbb C$, to a map $g$.
From lemma 2, $g$ is $J$-holomorphic.
$g$ is not constant while $|g^{'}(0)|=\displaystyle\lim_{n \to \infty}|\tilde g_n^{'}(0)|=1$.
So $g$ is a non trivial $J$-complex line.
\smallskip
\hskip 12 cm $\square$
\enddemo

\end